\documentclass{ws-adsaa}

\def\R{{\mathbb{R}}}
\def\Z{{\mathbb{Z}}}
\def\V{{\mathcal{V}}}
\def\F{{\mathcal{F}}}
\def\E{{\mathcal{E}}}

\def\sign{sgn}
\def\sinc{sinc}

\begin{document}

%
\catchline{}{}{}{}{}
%

\markboth{Authors' Names}
{Instructions for Typing Manuscripts (Paper's Title)}

\title{Continuous empirical wavelets systems}

\author{J\'er\^ome Gilles}
\address{Department of Mathematics \& Statistics, San Diego State University, 5500 Campanile Dr, San Diego, CA 92182, USA\\
\email{jgilles@sdsu.edu}
\http{https://jegilles.sdsu.edu/}}

\maketitle

\begin{history}
\received{Day Month Year}
\revised{Day Month Year}
\end{history}

\begin{abstract}
The recently proposed empirical wavelet transform was based on a particular type of filter. In this paper, we aim to propose a general framework for the construction of empirical wavelet systems in the continuous case. We define a well-suited formalism and then investigate some general properties of empirical wavelet systems. In particular, we provide some sufficient conditions to the existence of a reconstruction formula. In the second part of the paper, we propose the construction of empirical wavelet systems based on some classic mother wavelets.
\end{abstract}


\section{Introduction}
Originally, wavelets were developed to leverage some limitations of the short-term Fourier transform \cite{Daubechies1992,Mallat2009}. The primary goal was to develop families of functions corresponding to time-frequency atoms of different sizes, depending on their location in the time-frequency plane \cite{Daubechies1990}. This property permits to improve the accuracy of the extracted intrinsic characteristics of a function or signal. In the last decades, wavelets have proven to be very successful not only in mathematics (for instance to characterize certain function spaces \cite{Vedel2009,Triebel2,Triebel3}) but in almost all domain of science \cite{Jaffard2001} and especially in signal/image processing \cite{Startck1994,Livens1997,Chambolle1998,Chai2007,Gilles2008a,Ocak2009,Tang2009a,Shen2010,Hramov2015,Gilles2016a}. 
Despite the improvement provided by wavelets, the Gabor-Heinsenberg uncertainty principle \cite{Donoho1989} still limits the accuracy of the obtained time-frequency representation. In order to leverage this limitation and to extract a more detailed instantaneous frequency information, several alternative approaches were studied in the literature, we refer the reader to the review articles \cite{Boashash1992,Boashash1992a}. Among these new techniques, data-driven methods have recently received a lot of attention. In particular, the Hilbert-Huang Transform \cite{Huang1998} (HHT) has been widely studied. The authors aimed to take advantage of the property that the Hilbert transform easily permits to extract instantaneous amplitude and frequency of an Amplitude-Modulated/Frequency-Modulated signal (AM/FM). The HHT is then based on two steps: 1) the Empirical Mode Decomposition (EMD) is used to extract the intrinsic harmonic modes (the AM/FM components), 2) the Hilbert transform is applied to each mode to extract its corresponding instantaneous amplitudes and frequencies. Despite its success in a wide variety of applications, the HHT has a major drawback: the EMD step is a purely algorithmic method and lacks of solid mathematical background making its behavior difficult to predict and sensitive to noise. To overcome this drawback, a few alternatives inspired by the HHT were proposed in the literature. Hou et al. \cite{Hou2011} introduced a variational model based on $L^1$ minimization to directly extract the parameters of an AM/FM model. Dragomiretskiy et al. \cite{Dragomiretskiy2014} proposed to replace the EMD step by a method, called Variational Mode Decomposition (VMD), incorporating the notion of analytical function in a variational model to extract the harmonic modes of the signal. The synchrosqueezed wavelet transform \cite{Daubechies2011} performs a reallocation process on the standard continuous wavelet transform followed by a ridge detection step to build an accurate time-frequency representation. Another wavelet based alternative, called the Empirical Wavelet Transform (EWT), has been proposed in \cite{Gilles2013,Gilles2013a} to construct an adaptive, i.e data-driven, wavelet transform. The purpose of the EWT is to achieve an adaptive decomposition into harmonic modes of an input signal to replace the EMD step. The EWT, thanks to the well established wavelet theory, expects to provide solid theoretical foundations. The main idea at the core of the EWT is based on the work in \cite{Flandrin2004a} which experimentally shows that the EMD behaves like an adaptive filter bank, i.e filters based on data-driven supports in the frequency domain. This property can be easily explained since the Fourier transform of an AM/FM component is of ``almost'' compact support (i.e most of its energy is concentrated within a certain frequency range). The first step of the EWT consists in detecting the position of these supports, i.e segment the Fourier spectrum of the input signal, to define a partition of the frequency domain. Next, a Littlewood-Paley type filter is defined on each support creating a wavelet filter bank driven by the given input signal. Finally, this filter bank is applied to the input signal to extract the harmonic modes contained within. The Hilbert transform can then be individually applied to these components to extract an accurate time-frequency information. If the EWT has already proven its efficacy to analyze signals from different fields of science and engineering \cite{Huang2018,Huang2019,lfm,wind,ecg,speech,singlephase,seismic,EEGEWT}, it currently proposes only one type of wavelet and its theoretical aspects remain to be investigated. In this paper, we aim to establish a general framework/formalism to design different empirical wavelets and prove some global properties in the continuous framework. In Section~\ref{sec:notations}, we recall some notations and basic definitions which will be used throughout the paper. In Section~\ref{sec:partition}, we define the formalism to manipulate arbitrary partitions of the real line (i.e corresponding to the detected harmonic mode supports). Section~\ref{sec:EWS} defines Empirical Wavelet Systems and investigates some of their basic properties like the necessary and sufficient conditions to the existence of a reconstruction formula. Different families of empirical wavelets and some of their specific properties are investigated in Section~\ref{sec:families}. Finally, this work will be concluded in Section~\ref{sec:conclusion}.

\section{Generalities - Notations}\label{sec:notations}
In this paper, we will consider that all functions belong to $L^1(\R)\cap L^2(\R)$ equipped with its usual inner-product, defined by
$$\forall f,g\in L^1(\R)\cap L^2(\R)\quad;\quad \langle f,g\rangle=\int_\R f(t)\overline{g(t)}dt.$$
The Fourier transform of $f\in L^1(\R)\cap L^2(\R)$ and its inverse are defined by 
$$\F(f)(\xi)=\hat{f}(\xi)=\int_\R f(t)e^{-2\pi\imath\xi t}dt,$$ 
and $$f(t)=\F^{-1}(\hat{f})(t)=\check{\hat{f}}(t)=\int_\R \hat{f}(\xi)e^{2\pi\imath \xi t}d\xi,$$
respectively. Throughout the paper, the variable $\xi\in\R$ will denote the frequency. 
We will also use the following standard operators:
\begin{itemize}
 \item modulation $E_a$: $\forall f\in L^1(\R)\cap L^2(\R),E_af(t)=e^{2\pi\imath at}f(t)$,
 \item translation $T_a$: $\forall f\in L^1(\R)\cap L^2(\R),T_af(t)=f(t-a)$,
 \item scaling $D_a$: $\forall a>0,\forall f\in L^1(\R)\cap L^2(\R),D_af(t)=\frac{1}{\sqrt{a}}f\left(\frac{t}{a}\right)$,
\end{itemize}
it is easy to check that these operators have the following properties:
\begin{align}
 \F(E_af)=T_a\hat{f}\quad &;\quad \F^{-1}\left(E_a\hat{f}\right)=T_{-a}f\\
 \F(T_af)=E_{-a}\hat{f}\quad &;\quad \F^{-1}\left(T_a\hat{f}\right)=E_af\\
 \F(D_af)=D_{1/a}\hat{f}\quad &;\quad \F^{-1}\left(D_a\hat{f}\right)=D_{1/a}f
\end{align}

\section{Partitioning of the Fourier line}\label{sec:partition}
In this section, we define the formalism which will be used to describe the Fourier supports used in the construction of empirical wavelets. Let $\V=\{\nu_n\}_{n=n_m}^{n_M}$ be a set of boundary points, where $n_m,n_M\in\Z$ and $n_m<n_M$. For instance, the algorithm based on the scale-space theory and described in \cite{Gilles2014a} can be used to find such partition. We adopt the convention $\nu_0=0$ which implies that if $n<0$ then $\nu_n<0$ and if $n>0$ then $\nu_n>0$. In most practical cases $\nu_0=0$ is excluded (this will correspond to have a low pass-filter in the filter bank), and in such case we will denote the partition $\V^*=\{\nu_n\}_{n=n_m,n\neq 0}^{n_M}$. 
\begin{figure}[!t]
\centering
\includegraphics[width=\textwidth]{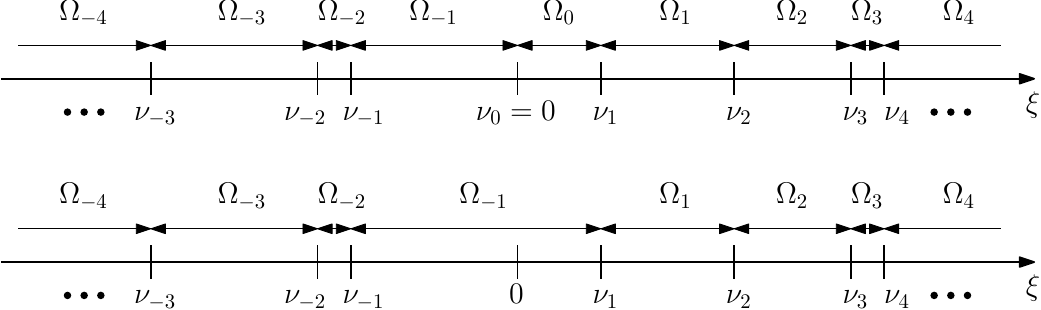}
\caption{Example of partitions $\V$ (top) and $\V^*$ (bottom) of the Fourier line.}
\label{fig:partition}
\end{figure}
A partition of the Fourier line based on the set of boundaries $\V$ can be defined as the set of intervals, or supports, of the form
$$\forall n=n_m,\ldots,n_M\quad,\quad \Omega_n=[\nu_n,\nu_{n+1}].$$
In the rest of the paper, we will denote $\Omega=\{\Omega_n\}_{n=n_m}^{n_M}$ a given partition of the Fourier line.
In the case of $\V^*$, the same definition is used except that $n=n_m,\ldots,n_M;n\neq 0$ and $\Omega_{-1}=[\nu_{-1},\nu_1]$ (i.e the zero frequency belongs to $\Omega_{-1}$).
An example for both cases is given in Figure~\ref{fig:partition}

\begin{figure}[!t]
\centering
\includegraphics[width=\textwidth]{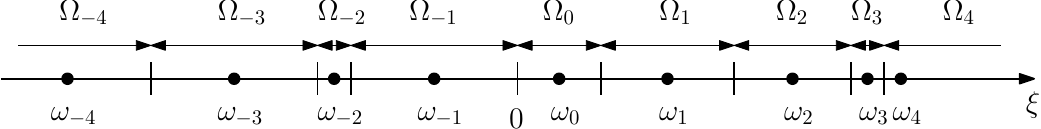}
\caption{Example of a $\V$ partition with infinite rays with its associated support centers.}
\label{fig:center}
\end{figure}

We can distinguish four main types of partitions (note that some of these cases are not exclusive):
\begin{enumerate}
 \item Infinite number of supports: $n_m=-\infty$ and $n_M=+\infty$ (i.e $n\in\Z$),
 \item Finite number of supports: $n_m$ and $n_M$ are finite,
 \item Right ray: $\nu_{n_M}=+\infty$, i.e the far right support is of the form $\Omega_{n_M-1}=[\nu_{n_M-1},+\infty)$,
 \item Left ray: $\nu_{n_m}=-\infty$, i.e the far left support is of the form $\Omega_{n_m}=(-\infty,\nu_{n_m+1}]$.
\end{enumerate}
The length of a support will be denoted $|\Omega_n|=\nu_{n+1}-\nu_n$ with the straightforward adaptation for $\V^*$. We will only consider cases avoiding supports of length 
either zero or infinite (except eventually potential right and/or left rays). Next, we need to define the center of a support. If the considered support $\Omega_n$ is compact then its center is defined by 
$$\omega_n=\frac{\nu_{n+1}+\nu_n}{2},$$
otherwise, in the case of rays $\Omega_{n_m}$ and $\Omega_{n_M-1}$, we use the size of the adjacent compact support:
\begin{equation}\label{eq:leftcen}
\omega_{n_m}=\nu_{n_m+1}-\frac{|\Omega_{n_m+1}|}{2}=\frac{3\nu_{n_m+1}-\nu_{n_m+2}}{2}, 
\end{equation}
and
\begin{equation}\label{eq:rightcen}
\omega_{n_M-1}=\nu_{n_M-1}+\frac{|\Omega_{n_M-2}|}{2}=\frac{3\nu_{n_M-1}-\nu_{n_M-2}}{2}. 
\end{equation}
An example is given is figure \ref{fig:center}. Straightforwardly the same definition can be used in the case of $\V^*$ except for $\omega_{-1}$ where 
$$\omega_{-1}=\frac{\nu_{-1}+\nu_1}{2}.$$

\section{Empirical Wavelet Systems}\label{sec:EWS}
An Empirical Wavelet System (EWS) is a set of filters whose supports in the Fourier domain correspond to the supports provided by the partition of the Fourier line as described in 
the previous section. An EWS is defined in the following way:
\begin{definition}
Let $\psi\in L^1( \R)\cap L^2(\R)$ be a function such that its Fourier transform is localized around the zero frequency. An Empirical Wavelet System, denoted $\{\psi_{b,n}\}_{b\in\R, n=n_m,\ldots,n_M}$, generated by $\psi$ is defined in the Fourier domain by, $\forall \xi\in\R$,
\begin{equation}\label{eq:EWSa}
\widehat{\psi}_{b,n}(\xi) = E_{-b}T_{\omega_n}D_{a_n}\widehat{\psi}\left(\xi\right). 
\end{equation}
or equivalently in the time domain by, $\forall t\in \R$,
\begin{equation}
\psi_{b,n}(t) = T_{b}E_{\omega_n}D_{1/a_n}\psi(t)
\end{equation}
where $\omega_n$ is the center of the support $\Omega_n$ as defined in the previous section and $a_n$ is a scaling factor whose choice depends on the used function $\psi$ and must 
be studied on a case by case basis.
\end{definition}
Let us interpret this definition. Since we assume that $\hat{\psi}$ is localized around the zero frequency, equation~\eqref{eq:EWSa} implies that $\widehat{\psi}_{b,n}$ is 
localized around the center $\omega_n$ of the support $\Omega_n$. The scaling factor, $a_n$, has to be chosen such that the support of $\widehat{\psi}_{b,n}$ has a width of about $|\Omega_n|$ to make sure $\widehat{\psi}_{b,n}$ is mostly localized on $\Omega_n$. 
The modulation operator corresponds to translating $\psi_{b,n}$ in the time domain. Therefore the family of functions $\{\psi_{b,n}\}_{b\in\R, 
n=n_m,\ldots,n_M}$ corresponds to a set of bandpass filters associated to the given partition $\Omega$.\\
For notational convenience, in the remaining of the paper we will denote $\widehat{\psi}_n=T_{\omega_n}D_{a_n}\widehat{\psi} 
$ (i.e $\widehat{\psi}_{b,n}=E_{-b}\widehat{\psi}_n$).

Given an Empirical Wavelet System, we can now give a general definition of the Continuous Empirical Wavelet Transform.
\begin{definition}\label{def:cewt}
The Continuous Empirical Wavelet Transform (CEWT), generated by $\{\widehat{\psi}_{b,n}\}_{b\in\R, n=n_m,\ldots,n_M}$, of a real or complex-valued function $f\in L^1(\R)\cap L^2(\R)$ is given 
by,
\begin{equation}
\E_{\psi}^f(b,n)= \langle\widehat{f},E_{-b}T_{\omega_n}D_{a_n}\widehat{\psi}\rangle= \langle f,T_bE_{\omega_n}D_{1/a_n}\psi\rangle. 
\end{equation}
\end{definition}
Like the classic wavelet transform, the CEWT can be rewritten as a filtering process. This is given by the following proposition ($\star$ will denote the convolution and we will 
denote $\psi^*(t)\equiv \psi(-t)$),

\begin{proposition}\label{prop-cewtconv}
The Continuous Empirical Wavelet Transform $\E_{\psi}^f(b,n)$ is equivalent to the convolution of $f$ with the function $\overline{\psi^*_n}(t)$, i.e
\begin{equation}
\E_{\psi}^f(b,n) = \left(f\star\overline{\psi^*_n}\right)(b)= \F^{-1}\left(\widehat{f}\cdot\overline{\widehat{\psi}_n}\right)(b).
\end{equation}
\end{proposition}
\begin{proof}
From definition~\ref{def:cewt} of the CEWT, given functions $f$ and $\psi$, we have,
\begin{align*}
\E_{\psi}^f(b,n)= \left\langle f,T_bE_{\omega_n}D_{1/a_n}\psi \right\rangle &= \int_\R f(t)\overline{T_bE_{\omega_n}D_{1/a_n}\psi(t)}dt \\
&=\int_\R f(t)\overline{T_b\psi_n(t)}dt\\ 
&=\int_\R f(t)\overline{\psi_n(t-b)}dt\\
&=\int_\R f(t)\overline{\psi^*_n(b-t)}dt\\
&=\left(f\star\overline{\psi^*_n}\right)(b).
\end{align*}
This proves the first part of the statement. Now, noticing that,
\begin{equation*}
\widehat{\overline{\psi^*_n}}= \int_\R \overline{\psi_n(-t)}e^{-2\pi \imath \xi t}dt
=\overline{\int_\R \psi_n(-t)e^{2\pi \imath \xi t}dt}
=\overline{\int_\R \psi_n(t)e^{-2\pi \imath \xi t}dt}
=\overline{\widehat{\psi}_n(\xi)}, 
\end{equation*}
we can rewrite the convolution obtained above as a pointwise multiplication in the Fourier domain,
\begin{align*}
\E_{\psi}^f(b,n)=\F^{-1}\left(\F\left(f(t)\star\overline{\psi^*_n(t)}\right)(b)\right)= \F^{-1}\left(\widehat{f}\cdot\widehat{\overline{\psi^*_n}}\right)(b)= 
\F^{-1}\left(\widehat{f}\cdot\overline{\widehat{\psi_n}}\right)(b).
\end{align*}
This provides the second part of the statement and ends the proof.
\end{proof}
The following proposition provides a sufficient condition to reconstruct $f$ from $\E_{\psi}^f$ by defining a dual set of empirical wavelets.
\begin{proposition}\label{prop:icewt}
Assume that $0<\sum_{n=n_m}^{n_M}\left|\widehat{\psi_n}(\xi)\right|^2<\infty\; a.e$, let us define the set of dual empirical wavelets $\{\phi_n\}_{n=n_m}^{n_M}$ by
$$\forall \xi\in\R\quad,\quad\widehat{\phi_n}(\xi)=\frac{\widehat{\psi}_n(\xi)}{\sum_{n=n_m}^{n_M}\left|\widehat{\psi}_n(\xi)\right|^2}$$
then 
\begin{equation}
f(t) = \sum_{n=n_m}^{n_M}\left(\E_{\psi}^f(\cdot,n)\star\phi_n\right)(t). 
\end{equation}
\end{proposition}
\begin{proof}
Using the Fourier transform and its inverse, we can write
\begin{align*}
\sum_{n=n_m}^{n_M}\left(\E_{\psi}^f(\cdot,n)\star\phi_n\right)(t)&=\F^{-1}\left(\F\left(\sum_{n=n_m}^{n_M}\left(\E_{\psi}^f(\cdot,n)\star\phi_n\right)(t)\right)\right)\\
&=\F^{-1}\left(\sum_{n=n_m}^{n_M}\F\left(\E_{\psi}^f(\cdot,n)\star\phi_n\right)(\xi)\right)\\
&=\F^{-1}\left(\sum_{n=n_m}^{n_M}\widehat{\E_{\psi}^f}(\xi,n)\widehat{\phi_n}(\xi)\right)\\
&=\F^{-1}\left(\sum_{n=n_m}^{n_M}\hat{f}(\xi)\overline{\widehat{\psi_n}}(\xi)\widehat{\phi_n}(\xi)\right)\\
&=\F^{-1}\left(\hat{f}(\xi)\sum_{n=n_m}^{n_M}\overline{\widehat{\psi_n}}(\xi)\frac{\widehat{\psi}_n(\xi)}{\sum_{n=n_m}^{n_M}\left|\widehat{\psi_n}(\xi)\right|^2}\right)\\
&=\F^{-1}\left(\hat{f}(\xi)\frac{\sum_{n=n_m}^{n_M}\left|\widehat{\psi}_n(\xi)\right|^2}{\sum_{n=n_m}^{n_M}\left|\widehat{\psi_n}(\xi)\right|^2}\right)\\
&=\F^{-1}\left(\hat{f}(\xi)\right)=f(t).
\end{align*}
This ends the proof.
\end{proof}
A particular case of the previous proposition is given by the following corollary.

\begin{corollary}\label{coro:tightewt}
If $\sum_{n=n_m}^{n_M}\left|\widehat{\psi_n}(\xi)\right|^2=A<\infty\; a.e$ then 
\begin{equation}
f(t) = \frac{1}{A}\sum_{n=n_m}^{n_M}\left(\E_{\psi}^f(\cdot,n)\star\psi_n\right)(t). 
\end{equation}
\end{corollary}
\begin{proof}
We apply Proposition~\ref{prop:icewt} with 
$$\forall \xi\in\R\quad,\quad \widehat{\phi_n}(\xi)=\frac{\widehat{\psi}_n(\xi)}{\sum_{n=n_m}^{n_M}\left|\widehat{\psi}_n(\xi)\right|^2}=\frac{\widehat{\psi}_n(\xi)}{A}$$
to get the expected result.
\end{proof}

\section{Construction of empirical wavelet systems}\label{sec:families}
In this section we study some families of empirical wavelets based on classic mother wavelets.
\subsection{Empirical Littlewood-Paley wavelets}
\begin{figure}[!t]
\centering\includegraphics[width=\textwidth]{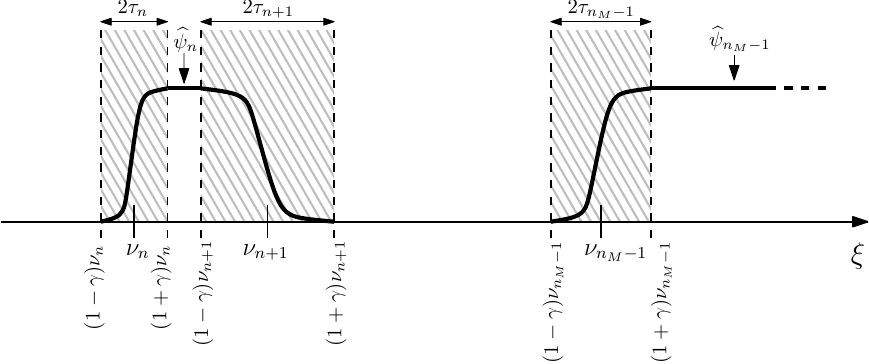}
\caption{Construction of empirical Littlewood-Paley wavelets}
\label{fig:meyer}
\end{figure}

\begin{figure}[!t]
\centering\includegraphics[width=\textwidth]{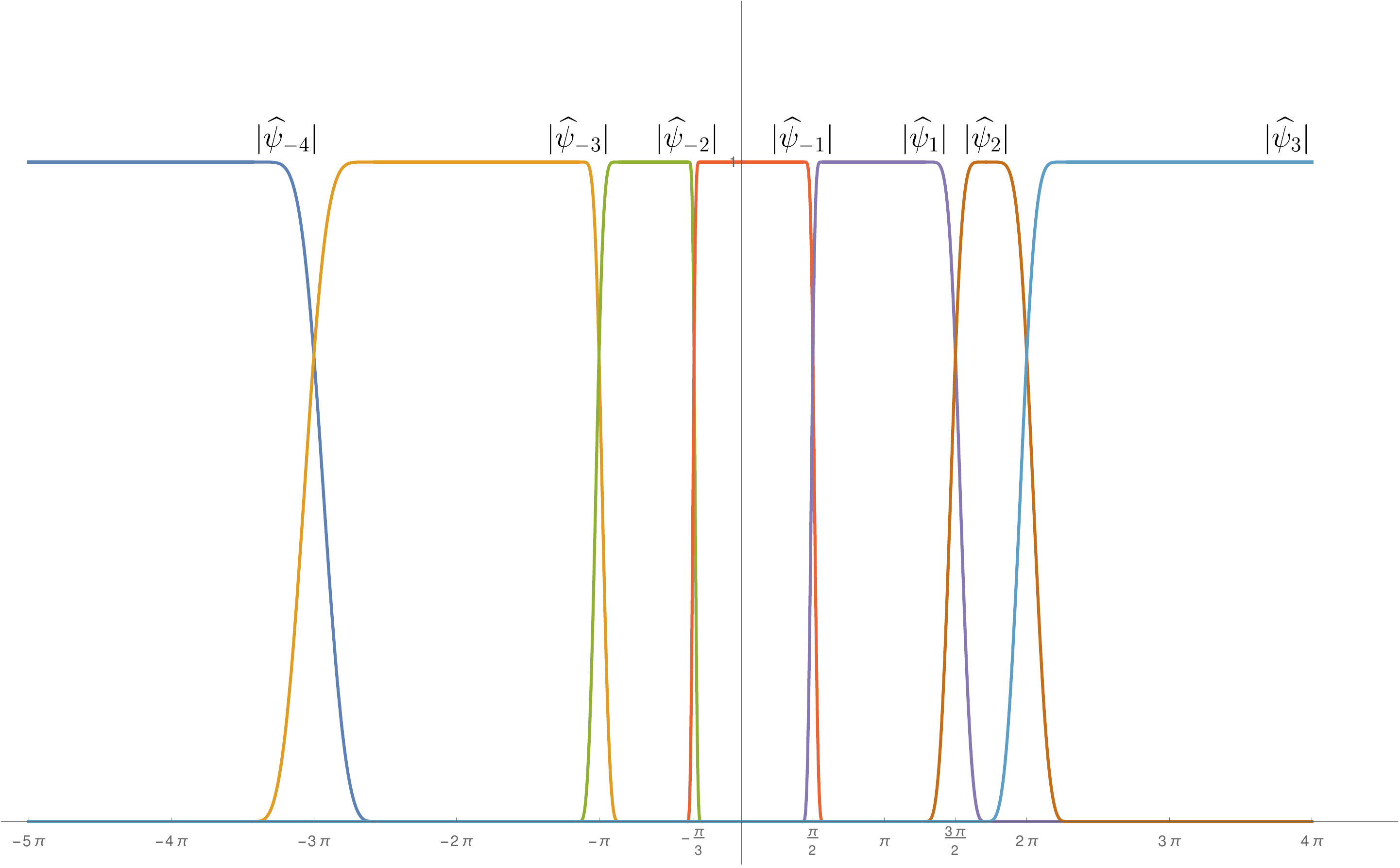}
\caption{Example of an empirical Littlewood-Paley wavelet filter bank in the Fourier domain corresponding to a $\V^*$ partition with left and right rays.}
\label{fig:exlp}
\end{figure}

Here, we revisit Littlewood-Paley (LP) wavelets proposed in the original EWT paper \cite{Gilles2013} but using the formalism 
defined in the previous sections. The definition we provide next is more general (as opposed to the original definition \cite{Gilles2013}) in the sense that it covers the case 
of complex signals, i.e non symmetric set of boundaries. The main idea in the construction of a LP-EWS is to define transition intervals centered at each boundary, this is illustrated in 
Figure~\ref{fig:meyer}. The width of the transition area centered at $\nu_n$ will be denoted $2\tau_n$ and in practice we will choose $\tau_n=\gamma\nu_n$ (the special case $n=0$ 
will be discussed hereafter) where $\gamma$ is a constant to be determined. Following the same argument as in the original work \cite{Gilles2013}, it can be shown that Corollary~\ref{coro:tightewt} can be fulfilled with $A=1$ by choosing $\gamma$ such that two consecutive transition intervals do not overlap.
First, let consider the case of a partition $\V$, i.e when $\nu_0=0$ is in the set of boundaries. On compact intervals $\Omega_n$, the empirical Littlewood-Paley wavelet is defined by
$$\widehat{\psi}_n^{LP}(\xi)=
\begin{cases}
  1 \qquad\qquad \quad\; \text{if}\quad (1+\sign(\nu_n)\gamma)\nu_n\leq\xi\leq(1-\sign(\nu_{n+1})\gamma)\nu_{n+1},\\
  \cos \left[\frac{\pi}{2}\beta \left(\frac{1}{2\gamma|\nu_{n+1}|}\left(\xi-(1-\sign(\nu_{n+1})\gamma)\nu_{n+1}\right)\right)\right]\\ 
  \qquad\qquad\; \text{if}\quad (1-\sign(\nu_{n+1})\gamma)\nu_{n+1}\leq \xi\leq (1+\sign(\nu_{n+1})\gamma)\nu_{n+1},\\
  \sin \left[\frac{\pi}{2}\beta \left(\frac{1}{2\gamma|\nu_n|}\left(\xi-(1-\sign(\nu_n)\gamma)\nu_n\right)\right)\right]\\ 
  \qquad\qquad\qquad\; \text{if}\quad (1-\sign(\nu_n)\gamma)\nu_n\leq \xi\leq (1+\sign(\nu_n)\gamma)\nu_n,\\
  0 \qquad\qquad\qquad\qquad \qquad\qquad \text{otherwise},
\end{cases}$$
where $\beta(x)=x^4(35-84x+70x^2-20x^3)$ (this is a classic choice used in the construction of Meyer wavelets \cite{Daubechies1992}). Littlewood-Paley wavelets on left and right rays are respectively 
defined by
$$\widehat{\psi}_{n_m}^{LPl}(\xi)=
\begin{cases}
  1 \qquad\qquad \qquad\qquad\qquad\qquad  \text{if}\quad \xi\leq(1-\sign(\nu_{n_m+1})\gamma)\nu_{n_m+1},\\
  \cos\left[\frac{\pi}{2}\beta\left(\frac{1}{2\gamma|\nu_{n_m+1}|}\left(\xi-(1-\sign(\nu_{n_m+1})\gamma)\nu_{n_m+1}\right)\right)\right]\\ 
  \quad \text{if}\quad (1-\sign(\nu_{n_m+1})\gamma)\nu_{n_m+1}\leq \xi\leq (1+\sign(\nu_{n_m+1})\gamma)\nu_{n_m+1},\\
  0 \qquad\qquad\qquad\qquad \qquad\qquad \text{otherwise},
\end{cases}$$
and
$$\widehat{\psi}_{n_M-1}^{LPr}(\xi)=
\begin{cases}
  1 \qquad\qquad \qquad\qquad\qquad  \text{if}\quad \xi\geq(1+\sign(\nu_{n_M-1})\gamma)\nu_{n_M-1},\\
  \sin\left[\frac{\pi}{2}\beta\left(\frac{1}{2\gamma|\nu_{n_M-1}|}\left(\xi-(1-\sign(\nu_{n_M-1})\gamma)\nu_{n_M-1}\right)\right)\right]\\ 
  \;\text{if}\quad (1-\sign(\nu_{n_M-1})\gamma)\nu_{n_M-1}\leq \xi\leq (1+\sign(\nu_{n_M-1})\gamma)\nu_{n_M-1},\\
  0 \qquad\qquad\qquad\qquad\qquad \text{otherwise},
\end{cases}$$
where $\gamma$ must be chosen such that 
\begin{equation}\label{eq:gamma}
\gamma<\min_{n=n_m,\ldots,n_M,|\Omega_n|<\infty}\left(\frac{|\Omega_n|}{2|\omega_n|}\right),
\end{equation}
to avoid overlap of transition intervals. It remains to choose $\tau_0$ (since $\gamma\nu_0=0$), we propose to choose $\tau_0=\min(\tau_{-1},\tau_1)$ to fulfill the conditions.\\
Second, we consider the case of a partition $\V^*$, i.e $\nu_0=0$ is not used in the definition of the partition. The same definitions as for the $\V$ case can be reused with of course 
the straightforward adaptation for $\widehat{\psi}_{-1}^{LP}$ which must use the boundaries $\nu_{-1}$ and $\nu_1$. However, a special case can occur when $\omega_{-1}=0$ (i.e 
$\nu_{-1}=-\nu_1$) which will be an issue in equation \eqref{eq:gamma}. In order to avoid this problem, we propose to choose $\gamma$ such that
\begin{equation}\label{eq:gammabis}
\gamma<\min\left(\min_{n=n_m,\ldots,n_M,n\neq-1,|\Omega_n|<\infty}\left(\frac{|\Omega_n|}{2|\omega_n|}\right),\frac{1}{2}\right),
\end{equation}
($1/2$ corresponds to consider half the width of $\Omega_{-1}$). Following the same proof as in Proposition 1 in the article \cite{Gilles2013}, it is easy to check that a 
family of empirical Littlewood-Paley wavelets fulfills Corollary~\ref{coro:tightewt} with bound $A=1$.\\
An example of an empirical Littlewood-Paley wavelets filter bank is given in Figure~\ref{fig:exlp} for the set of 
boundaries $\V^*=\{-\infty,-3\pi,-\pi,-\pi/3,\pi/2,3\pi/2,2\pi\}$.

\subsection{Empirical Meyer wavelets}
\begin{figure}[!t]
\centering\includegraphics[width=\textwidth]{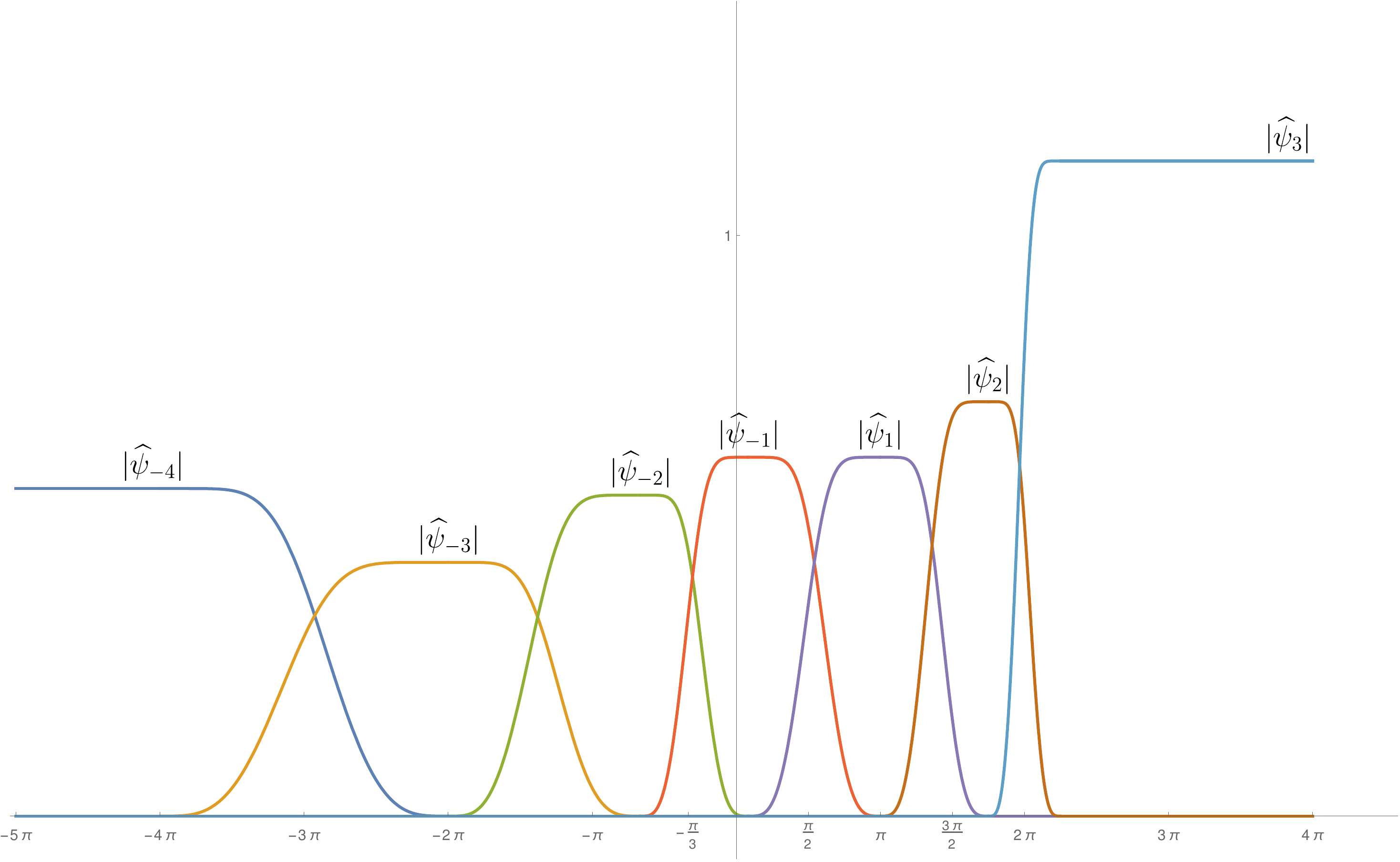}
\caption{Example of an empirical Meyer wavelet filter bank in the Fourier domain corresponding to a $\V^*$ partition with left and right rays.}
\label{fig:exmeyer}
\end{figure}

Meyer wavelets were an important example since this construction shows that orthogonal wavelets exist (while Meyer was trying to prove the opposite!). It is important to notice that the orthogonality is obtained thanks to a ``magic'' cancellation trick of phases of consecutive wavelets. Unfortunately, this trick does no longer work with the construction of empirical Meyer wavelets given hereafter making them non-orthogonal. 
In the empirical case, we cannot directly write the empirical Meyer wavelet in the general form $\hat{\psi}_n=T_{\omega_n}D_{a_n}\hat\psi$ since $\psi$ is defined piecewise and we need to control the scaling factor independently for each sub-interval. Therefore, we provide a direct formulation for $\hat{\psi}_n$:
$$\hat{\psi}_n^M(\xi)=\sqrt{\frac{2}{\omega_{n+1}-\omega_{n-1}}}e^{\imath\frac{4\pi}{3\max(|\omega_{n-1}|,|\omega_{n+1}|)}}
\begin{cases} 
\sin\left(\frac{\pi}{2}\beta\left(\frac{\xi-\omega_{n-1}}{\omega_n-\omega_{n-1}}\right)\right)\\  
\qquad\qquad\qquad \text{if}\; \omega_{n-1}\leq \xi\leq \omega_n \\ 
\cos\left(\frac{\pi}{2}\beta\left(\frac{\xi-\omega_n}{\omega_{n+1}-\omega_n}\right)\right)\\
\qquad\qquad\qquad \text{if}\; \omega_n\leq \xi\leq \omega_{n+1} 
\end{cases}.
$$
The rays cases are given by
$$\hat{\psi}_{n_m}^{Ml}(\xi)=\sqrt{\frac{2}{\omega_{{n_m}+1}-\omega_{n_m}}}e^{\imath\frac{4\pi}{3|\omega_{{n_m}+1}|}}
\begin{cases} 
1 \qquad\qquad\qquad \text{if}\; \xi\leq\omega_{n_m}\\
\cos\left(\frac{\pi}{2}\beta\left(\frac{\xi-\omega_{n_m}}{\omega_{n_m+1}-\omega_{n_m}}\right)\right)\\
\qquad\qquad\qquad\;\, \text{if}\; \omega_{n_m}\leq \xi\leq \omega_{n_m+1} \\
0 \qquad\qquad\qquad \text{if}\; \omega_{n_m+1}\leq \xi
\end{cases}.
$$
and
$$\hat{\psi}_{n_M-1}^{Mr}(\xi)=\sqrt{\frac{2}{\omega_{n_M-1}-\omega_{n_M-2}}}e^{\imath\frac{4\pi}{3|\omega_{n_M-2}|}}
\begin{cases} 
0 \qquad\qquad \text{if}\; \xi\leq\omega_{n_M-2}\\
\sin\left(\frac{\pi}{2}\beta\left(\frac{\xi-\omega_{n_M-2}}{\omega_{n_M-1}-\omega_{n_M-2}}\right)\right)\\
\qquad\qquad\; \text{if}\; \omega_{n_M-2}\leq \xi\leq \omega_{n_M-1} \\
1 \qquad\qquad \text{if}\; \omega_{n_M-1}\leq \xi
\end{cases}.
$$
An example of an empirical Meyer wavelet filter bank is given in Figure~\ref{fig:exmeyer}.
Empirical Meyer wavelets have the property that each wavelet is normalized.
\begin{proposition}
Each empirical Meyer wavelet (except the rays ones) is normalized, i.e $\forall n=n_m+1,\ldots,n_M-1,\|\psi_n^M\|_{L^2}=1$.
\end{proposition}
\begin{proof}
We have
\begin{align*}
\|\psi_n^M\|_{L^2}=\|\widehat{\psi}_n^M\|_{L^2}&=\frac{2}{\omega_{n+1}-\omega_{n-1}}\left(\int_{\omega_{n-1}}^{\omega_n}\sin^2\left[\frac{\pi}{2}\beta\left(\frac{\xi-\omega_{n-1}}{\omega_n-\omega_{n-1}}\right)\right]d\xi\right.\\ 
& \qquad\qquad\qquad\qquad\left.+\int_{\omega_n}^{\omega_{n+1}}\cos^2\left[\frac{\pi}{2}\beta\left(\frac{\xi-\omega_n}{\omega_{n-1}-\omega_n}\right)\right]d\xi\right)\\
&=\frac{2}{\omega_{n+1}-\omega_{n-1}}\left(\int_0^1\sin^2\left[\frac{\pi}{2}\beta(u)\right](\omega_n-\omega_{n-1})du\right.\\
& \qquad\qquad\qquad\qquad\left.+\int_0^1\cos^2\left[\frac{\pi}{2}\beta(u)\right](\omega_{n+1}-\omega_n)du\right)\\
&=\frac{2}{\omega_{n+1}-\omega_{n-1}}\left(\int_0^1(\omega_n-\omega_{n-1})\sin^2\left[\frac{\pi}{2}\beta(u)\right]du\right.\\
& \quad\qquad\qquad+\int_0^1(\omega_n-\omega_{n-1})\cos^2\left[\frac{\pi}{2}\beta(u)\right]du\\
& \quad\qquad\qquad\left.+\int_0^1(\omega_{n+1}-2\omega_n+\omega_{n-1})\cos^2\left[\frac{\pi}{2}\beta(u)\right]du\right)\\
&=\frac{2}{\omega_{n+1}-\omega_{n-1}}\bigg(\omega_n-\omega_{n-1}\\
& \quad\qquad\qquad\left.+(\omega_{n+1}-2\omega_n+\omega_{n-1})\int_0^1\cos^2\left[\frac{\pi}{2}\beta(u)\right]du\right),
\end{align*}
as shown \cite{Daubechies1992}, $\int_0^1\cos^2\left[\frac{\pi}{2}\beta(u)\right]du=\frac{1}{2}$ hence
$$\|\psi_n^M\|_{L^2}=\frac{2}{\omega_{n+1}-\omega_{n-1}}\frac{2\omega_n-2\omega_{n-1}+\omega_{n+1}-2\omega_n+\omega_{n-1}}{2}=1.$$
\end{proof}

\begin{proposition}
The empirical Meyer wavelet system fulfill the following condition:
$$A\leq \sum_{n=n_m}^{n_M-1}|\hat{\psi}(\xi)|^2 \leq B$$
where
\begin{align*}
A&=\min\left(\frac{2}{\omega_{n_M-1}-\omega_{n_M-2}};\frac{2}{\omega_{n_m+1}-\omega_{n_m}};\min_{n=n_m+1,\ldots,n_M-2}\left(\frac{2}{\omega_{n+1}-\omega_{n-1}}\right)\right),
\end{align*}
and
\begin{align*}
B&=\max\left(\frac{2}{\omega_{n_M-1}-\omega_{n_M-2}};\frac{2}{\omega_{n_m+1}-\omega_{n_m}};\max_{n=n_m+1,\ldots,n_M-2}\left(\frac{2}{\omega_{n+1}-\omega_{n-1}}\right)\right). 
\end{align*}
\end{proposition}
\begin{proof}
By construction, we know that in the Fourier domain, only two consecutive filters have overlapping supports. It is straightforward from the definition of the rays to see that
$$\sum_{n=n_m}^{n_M-1}|\hat{\psi}(\xi)|^2=\begin{cases}
                                           \frac{2}{\omega_{n_m+1}-\omega_{n_m}}\qquad&\text{if}\; \xi\leq \omega_{n_m}\\
                                           \frac{2}{\omega_{n_M-1}-\omega_{n_M-2}} &\text{if}\; \xi\geq \omega_{n_M-1}
                                          \end{cases}.$$
Next, $\forall n=n_m+1,\ldots,n_M-2,\forall \xi\in [\omega_n,\omega_{n+1}]$, we have
\begin{align*}
\sum_{n=n_m}^{n_M-1}|\hat{\psi}(\xi)|^2&=|\hat{\psi}_n(\xi)|^2+|\hat{\psi}_{n+1}(\xi)|^2\\
&=\frac{2}{\omega_{n+1}-\omega_{n-1}}\cos^2\left[\frac{\pi}{2}\beta\left(\frac{\xi-\omega_n}{\omega_{n+1}-\omega_n}\right)\right]\\
&\qquad\qquad+\frac{2}{\omega_{n+2}-\omega_{n}}\sin^2\left[\frac{\pi}{2}\beta\left(\frac{\xi-\omega_n}{\omega_{n+1}-\omega_n}\right)\right].\\
\end{align*}
To simplify the notations, we will denote $u=\frac{\pi}{2}\beta\left(\frac{\xi-\omega_n}{\omega_{n+1}-\omega_n}\right)$. Thus, we have
\begin{align*}
\sum_{n=n_m}^{n_M-1}|\hat{\psi}(\xi)|^2&=\frac{2}{\omega_{n+1}-\omega_{n-1}}\cos^2(u)+\frac{2}{\omega_{n+2}-\omega_{n}}\sin^2(u)\\
&=\frac{2}{\omega_{n+1}-\omega_{n-1}}+\underbrace{\left(\frac{2}{\omega_{n+2}-\omega_{n}}-\frac{2}{\omega_{n+1}-\omega_{n-1}}\right)}_{\Theta}\sin^2(u).
\end{align*}
Given that $0\leq\sin^2(u)\leq 1$, we consider the different cases:
\begin{itemize}
 \item $\Theta>0$, we get
 $$\frac{2}{\omega_{n+1}-\omega_{n-1}}\leq \sum_{n=n_m}^{n_M-1}|\hat{\psi}(\xi)|^2\leq \frac{2}{\omega_{n+2}-\omega_{n}},$$
 \item $\Theta=0$, we get
 $$\sum_{n=n_m}^{n_M-1}|\hat{\psi}(\xi)|^2= \frac{2}{\omega_{n+1}-\omega_{n-1}},$$
 \item $\Theta<0$, we get
 $$\frac{2}{\omega_{n+2}-\omega_{n}}\leq \sum_{n=n_m}^{n_M-1}|\hat{\psi}(\xi)|^2\leq \frac{2}{\omega_{n+1}-\omega_{n-1}}.$$
\end{itemize}
Therefore,
\begin{align*}
\min\left(\frac{2}{\omega_{n+1}-\omega_{n-1}};\frac{2}{\omega_{n+2}-\omega_{n}}\right)\leq &\sum_{n=n_m}^{n_M-1}|\hat{\psi}(\xi)|^2 \\
&\leq \max\left(\frac{2}{\omega_{n+1}-\omega_{n-1}};\frac{2}{\omega_{n+2}-\omega_{n}}\right) 
\end{align*}
It remains to take the minimum and maximum over all possible intervals to get 
$$A\leq \sum_{n=n_m}^{n_M-1}|\hat{\psi}(\xi)|^2\leq B,$$
where
\begin{align*}
A&=\min\left(\frac{2}{\omega_{n_M-1}-\omega_{n_M-2}};\frac{2}{\omega_{n_m+1}-\omega_{n_m}};\min_{n=n_m+1,\ldots,n_M-2}\left(\frac{2}{\omega_{n+1}-\omega_{n-1}}\right)\right),
\end{align*}
and
\begin{align*}
B&=\max\left(\frac{2}{\omega_{n_M-1}-\omega_{n_M-2}};\frac{2}{\omega_{n_m+1}-\omega_{n_m}};\max_{n=n_m+1,\ldots,n_M-2}\left(\frac{2}{\omega_{n+1}-\omega_{n-1}}\right)\right). 
\end{align*}
This conclude the proof.
\end{proof}
We can notice that $A>0$ and $B<\infty$ which, thanks to Proposition~\ref{prop:icewt}, guarantee the existence of a dual empirical wavelet system and hence a reconstruction formula.

\subsection{Empirical Shannon wavelets}
\begin{figure}[!t]
\centering\includegraphics[width=\textwidth]{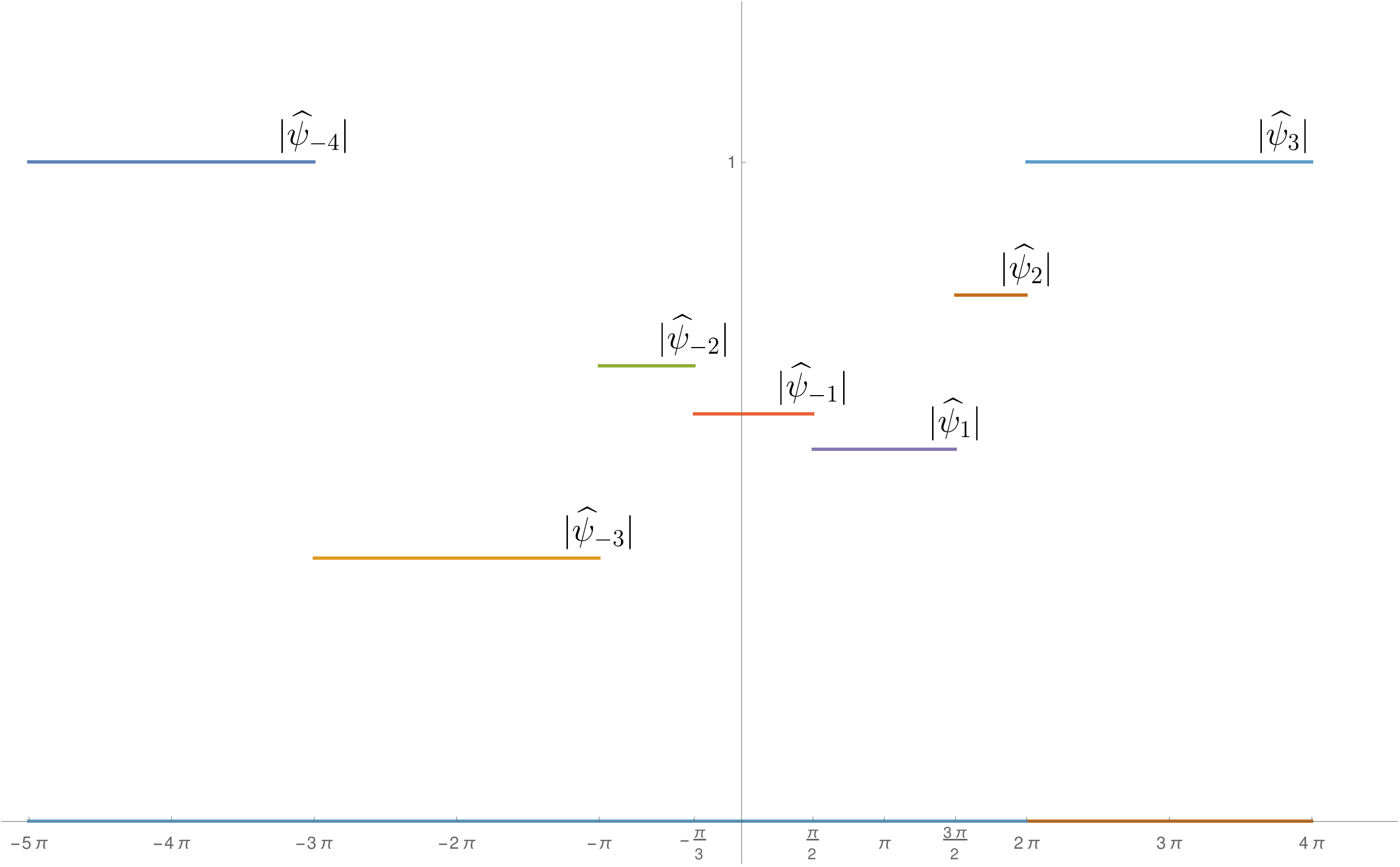}
\caption{Example of an empirical Shannon wavelet filter bank in the Fourier domain corresponding to a $\V^*$ partition with left and right rays.}
\label{fig:shannon}
\end{figure}

Shannon wavelets are commonly used in the Shannon sampling theory \cite{Shannon1}. These wavelets are made of $\sinc$ functions which correspond to characteristic functions over intervals in the 
Fourier domain. The construction of empirical Shannon wavelets is then very easy since the mother wavelet is defined by $$\widehat{\psi}^{SH}(\xi)=e^{-\imath\frac{\pi}{2}\left(\xi+\frac{3}{2}\right)}\chi_{[-1/2,1/2]}(\xi)=\begin{cases} e^{-\imath\xi/2}\quad &\text{if}\;\xi\in [-1/2,1/2)\\ 0 &\text{otherwise}\end{cases}.$$ 
Therefore, in order to keep the width of the scaled wavelet to be $|\Omega_n|$, we have to set $a_n=|\Omega_n|$ and thus the set of empirical Shannon wavelets is defined by
$$\widehat{\psi}_n^{SH}(\xi)=T_{\omega_n}D_{|\Omega_n|}\widehat{\psi}^{SH}(\xi)=\frac{1}{\sqrt{|\Omega_n|}}e^{-\imath\frac{\pi}{2}\left(\frac{\xi-\omega_n}{|\Omega_n|}+\frac{3}{2}\right)}\chi_{\left[\omega_n-\frac{|\Omega_n|}{2};\omega_n+\frac{|\Omega_n|}{2}\right)}(\xi).$$
Rays are simply defined by
$$\widehat{\psi}_n^{SHl}(\xi)=e^{-\imath\pi}\chi_{\left(-\infty,\nu_{n_m+1}\right)}(\xi),$$
and
$$\widehat{\psi}_n^{SHr}(\xi)=e^{-\imath\frac{\pi}{2}}\chi_{\left[\nu_{n_M-1},+\infty\right)}(\xi).$$
Examples of empirical Shannon filters for the same partition we used before is given in Figure~\ref{fig:shannon}
Empirical Shannon wavelet can be characterized by
\begin{proposition}
A set of empirical Shannon wavelets $\{\psi^{SH}\}$ fulfills the condition $\forall \xi\in\R$ (taking the convention that $|\Omega_n|=1$ for eventual left or right rays),
$$\frac{1}{\max_{n=n_m,\ldots,n_M} |\Omega_n|}\leq \sum_{n=n_m}^{n_M}\left|\widehat{\psi}_n^{SH}(\xi)\right|^2\leq \frac{1}{\min_{n=n_m,\ldots,n_M} |\Omega_n|}$$
\end{proposition}
\begin{proof}
By construction, the Shannon wavelets filters do not overlap in the Fourier domain hence $\forall n=n_m,\ldots,n_M,\forall \xi\in\Omega_n,\left|\widehat{\psi}_n^{SH}(\xi)\right|^2=\frac{1}{|\Omega_n|}$ or $\left|\widehat{\psi}_n^{SH}(\xi)\right|^2=1$ for left or right rays. Therefore the lower and upper bounds across all $n$ are given by the bounds provided in the proposition.
\end{proof}

\subsection{Empirical Gabor wavelets}
\begin{figure}[!t]
\centering\includegraphics[width=\textwidth]{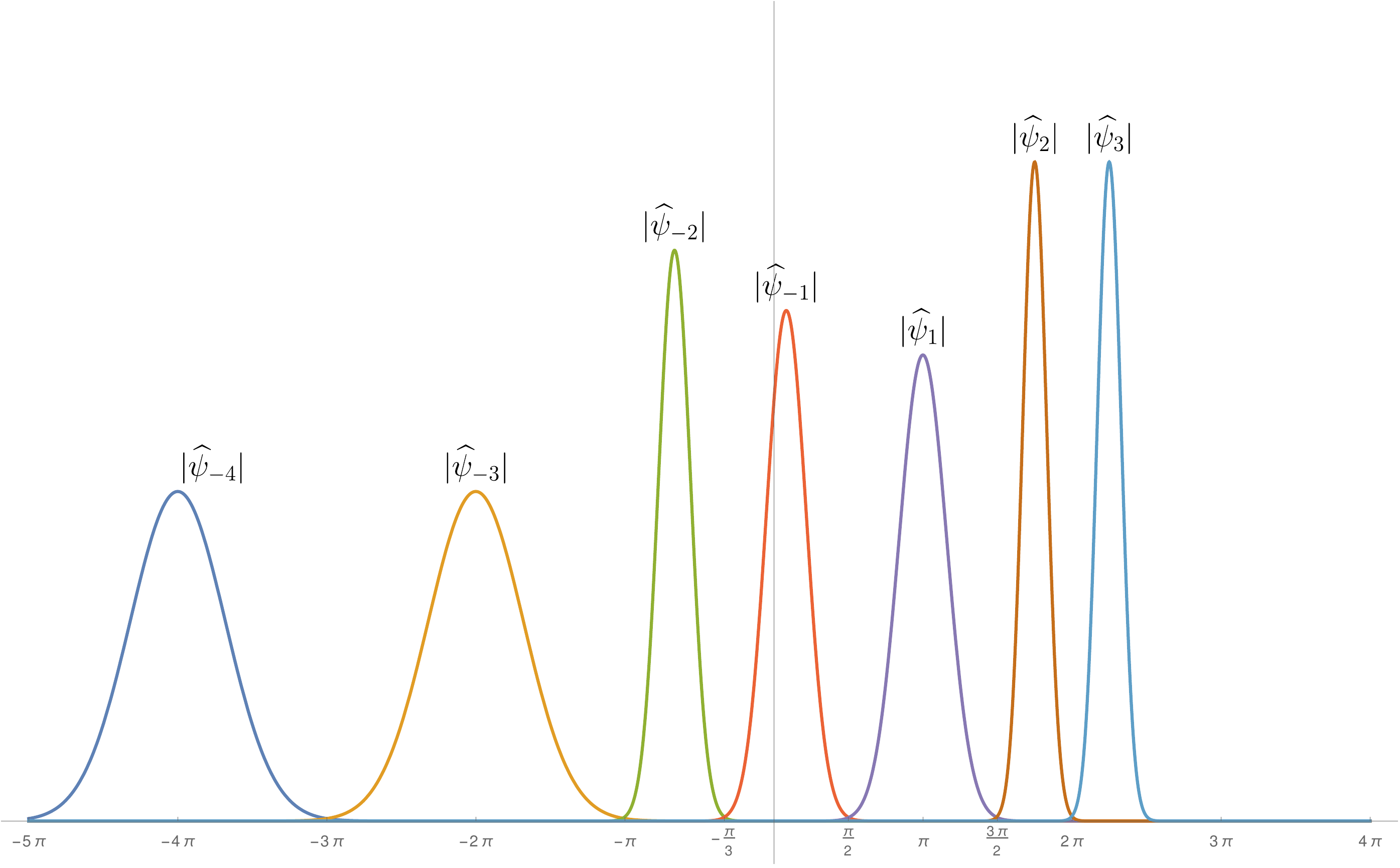}
\caption{Example of an empirical Gabor wavelets filter bank in the Fourier domain where option 1 was used for the rays.}
\label{fig:gabor1}
\end{figure}

\begin{figure}[!t]
\centering\includegraphics[width=\textwidth]{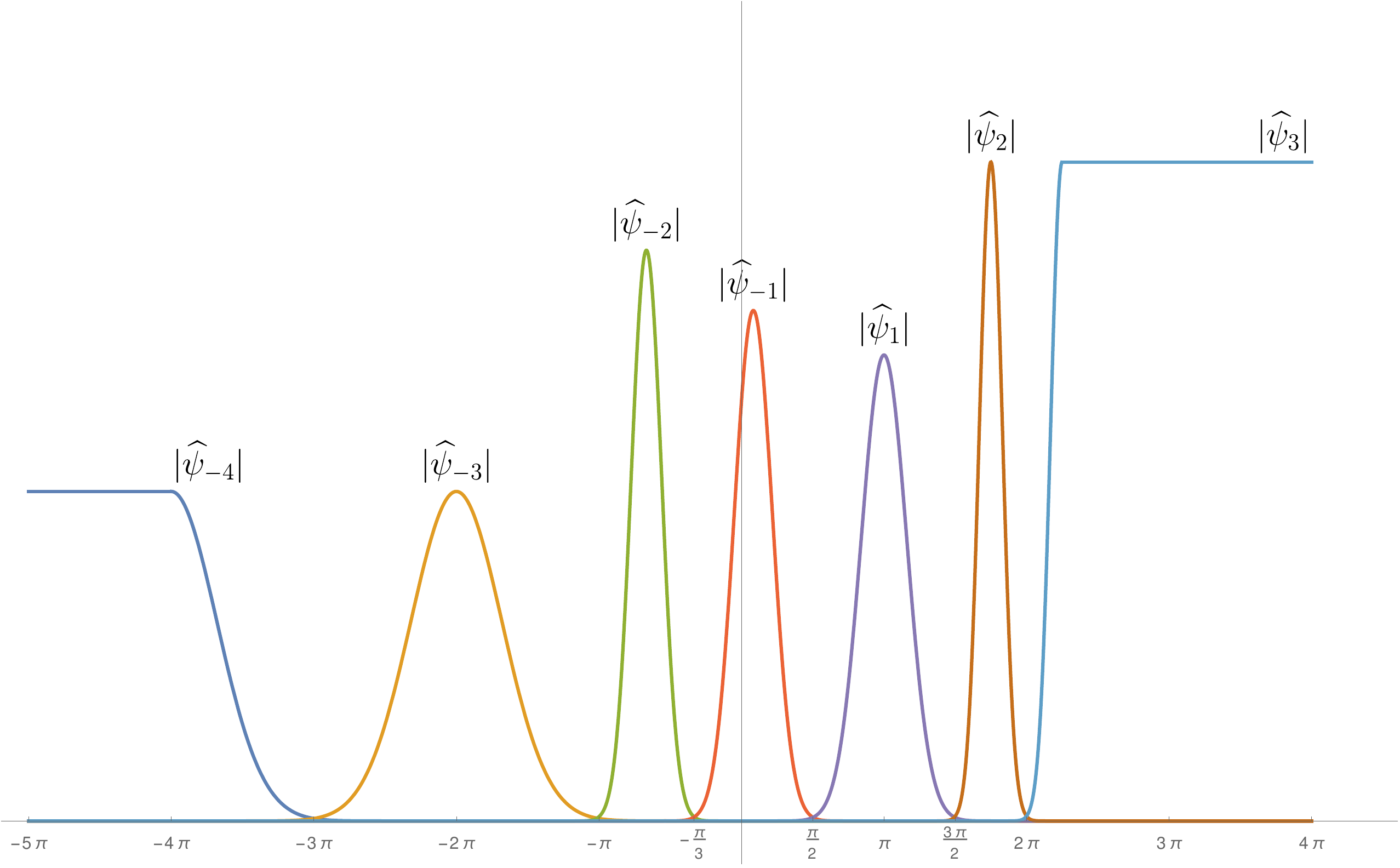}
\caption{Example of an empirical Gabor wavelets filter bank in the Fourier domain where option 2 was used for the rays.}
\label{fig:gabor2}
\end{figure}

The previous cases used mother wavelets having compact supports in the Fourier domain making straightforward the construction of empirical families. In this section, we address 
the case of Gabor wavelets \cite{Christensen2001,Christensen2010,Balazs2011,Ricaud2013} which do not have such property. We use the mother wavelet $\psi^G$ defined by $\forall\xi\in\R,\widehat{\psi}^G(\xi)=e^{-\pi(2.5\xi)^2}$, this choice 
ensure that $\widehat{\psi}^G$ is mostly localized in the interval $[-1/2,1/2]$ (in the sense that $99.999\%$ of $\|\psi^G\|_{L^2}$ comes from this interval) and that 
$\widehat{\psi}^G(0)=1$. In order to keep $\widehat{\psi}_n^G$ localized in $\Omega_n$, we need to choose $a_n=|\Omega_n|$. Therefore we define 
$$\forall \xi\in\R\quad,\quad \widehat{\psi}_n^G(\xi)=\frac{1}{\sqrt{|\Omega_n|}}e^{-\pi\left(\frac{2.5(\xi-\omega_n)}{|\Omega_n|}\right)^2},$$
and $\widehat{\psi}_{b,n}^G(\xi)=e^{-2\imath b\xi}\widehat{\psi}_n^G(\xi)$.
If rays are part of the partition, we propose two options to define them. The first option consists in defining one Gaussian for each ray in the following way:
\begin{equation}
\widehat{\psi}_{n_m}^{Gl}(\xi)=\frac{1}{\sqrt{|\Omega_{n_m+1}|}}e^{-\pi\left(\frac{2.5(\xi-\omega_{n_m})}{|\Omega_{n_m+1}|}\right)^2},
\end{equation}
and
\begin{equation}
\widehat{\psi}_{n_M-1}^{Gr}(\xi)=\frac{1}{\sqrt{|\Omega_{n_M-2}|}}e^{-\pi\left(\frac{2.5(\xi-\omega_{n_M-1})}{|\Omega_{n_M-2}|}\right)^2},
\end{equation}
where $\omega_{n_m}$ and $\omega_{n_M-1}$ are given by \eqref{eq:leftcen} and \eqref{eq:rightcen}, respectively. This option has the drawback that the spectral information ``far'' on the left or right (i.e outside the intervals where the left and right Gaussian are localized) are completely wiped out from the analysis. This could potentially be an issue for some applications hence we propose a second option to define the rays filters:
\begin{equation}
\widehat{\psi}_{n_m}^{\widetilde{Gl}}(\xi)=
\begin{cases}
\widehat{\psi}_{n_m}^{Gl}(\xi) \qquad&\text{if}\; \omega_{n_m}\leq\xi\leq \nu_{n_m+1}\\
\frac{1}{\sqrt{|\Omega_{n_m+1}|}} &\text{otherwise},
\end{cases}
\end{equation}
and
\begin{equation}
\widehat{\psi}_{n_M-1}^{\widetilde{Gr}}(\xi)=
\begin{cases}
\widehat{\psi}_{n_M-1}^{Gr}(\xi) \qquad&\text{if}\; \nu_{n_M-1}\leq\xi\leq \omega_{n_M-1}\\
\frac{1}{\sqrt{|\Omega_{n_M-2}|}} &\text{otherwise}.
\end{cases}
\end{equation}
This option does not discard the far spectral information. Examples of empirical Gabor filters for each options for the same partition we used before are given in Figures~\ref{fig:gabor1} and \ref{fig:gabor2}, respectively.
We can prove the following property.
\begin{proposition}
All types of empirical Gabor wavelet sets fulfill the condition
$$\forall\xi\in\R\quad,\quad 0<\sum_{n=n_m}^{n_M}\left|\widehat{\psi}_n(\xi)\right|^2<\infty.$$
Moreover, if the empirical Gabor wavelets set is build either based on an infinite set of supports or a finite set of supports using the second option described above then
$$\forall \xi\in\R\quad,\quad \sum_{n=n_m}^{n_M}\left|\widehat{\psi}_n(\xi)\right|^2\geq e^{-\frac{25\pi}{8}}.$$
\end{proposition}
\begin{proof}
Let an arbitrary set of empirical Gabor wavelets $\{\psi_n\}$. By construction, it is obvious that $\sum_{n=n_m}^{n_M}\left|\widehat{\psi}_n(\xi)\right|^2>0$. To check the upper bound, let consider first the case of an infinite set of supports, we write
\begin{align*}
 \sum_{n=n_m}^{n_M}\left|\widehat{\psi}_n(\xi)\right|^2=\sum_{n=n_m}^{n_M}\left|T_{\omega_n}D_{a_n}\widehat{\psi}(\xi)\right|^2&=\sum_{n=n_m}^{n_M}\left|T_{\omega_n}D_{a_n}e^{-\pi(2.5\xi)^2}\right|^2\\
 &=\sum_{n=n_m}^{n_M}\frac{1}{a_n}e^{-2\pi\left(2.5\frac{\xi-\omega_n}{a_n}\right)^2}.
\end{align*}
Let denote $a_{min}=\min_{n=n_m,\ldots,n_M}a_n$ and $a_{max}=\max_{n=n_m,\ldots,n_M}a_n$. Noticing that, by construction, $a_{min}$ and $a_{max}$ are both finite and positive, we have 
\begin{align*}
\sum_{n=n_m}^{n_M}\left|\widehat{\psi}_n(\xi)\right|^2\leq \frac{1}{a_{min}}\sum_{n=n_m}^{n_M}e^{-2\pi\left(2.5\frac{\xi-\omega_n}{a_{max}}\right)^2}\leq\frac{1}{a_{min}}\int_\R e^{-2\pi\left(2.5\frac{\xi-\omega}{a_{max}}\right)^2}d\omega<\infty,
\end{align*}
therefore
$$\forall\xi\in\R\quad,\quad 0<\sum_{n=n_m}^{n_M}\left|\widehat{\psi}_n(\xi)\right|^2<\infty.$$
Now assume that the used partition is either infinite or finite and in this case rays are defined via the second option describes above. Then we have that $\forall\xi\in\R,\exists n_*\in \{n_m,\ldots,n_M\}$ such that $\xi\in\Omega_{n_*}$. On this interval, again by construction, we necessarily have 
$$\sum_{n=n_m}^{n_M}\left|\widehat{\psi}_n(\xi)\right|^2\geq \left|\widehat{\psi}_{n_*}(\xi)\right|^2\geq \left|\widehat{\psi}_{n_*}(\nu_{n_*})\right|^2=\left|\widehat{\psi}(1/2)\right|^2=e^{-\frac{25\pi}{8}}.$$
\end{proof}
Note that in the last part of this proof, a finite partition using option 1 to build the wavelet does not guarantee such constant lower bound. Indeed, on left or right rays, as $\xi$ goes to infinity $\widehat{\psi}_n\to 0$. This proposition is useful because it guaranties that, using Proposition~\ref{prop:icewt}, we can build a set of function $\{\phi_n\}$ to reconstruct the original function from the empirical wavelet transform information.

\section{Numerical implementation}
In this section, we consider discrete signals made of $N$ samples (we don't consider signals of infinite length here). Given that the Fourier transform is periodic, we only consider the Fourier domain within one period. For notation convenience, we use the normalized interval of frequencies $(-\pi,\pi]$, i.e $\V\subset (-\pi,\pi)$. All results and constructions of the different families of wavelets can be easily reformulated in the discrete case by sampling the variable $\xi$. The corresponding numerical implementation of the Littlewood-Paley, Meyer, Gabor and Shannon empirical wavelets systems are made available through the Empirical Wavelet Toolbox available at \url{https://www.mathworks.com/matlabcentral/fileexchange/42141-empirical-wavelet-transforms}.

\section{Conclusion}\label{sec:conclusion}
In this paper, we defined a general framework to build continuous empirical wavelet systems. We provided some conditions to guarantee well-behaved dual systems useful to reconstruct a function from its transform. Finally, we showed the construction of several empirical wavelet systems based on some classic mother wavelets, their numerical implementation is freely available in the Empirical Wavelet Toolbox. In terms of future investigations, we expect to derive conditions to build empirical wavelet frames as well as provide some general closed form for the corresponding frame bounds. Extensions to higher dimension are also of interest but will raise some important questions related to the geometrical properties of the detected supports in the Fourier domain and how to derive conditions for the existence of a reconstruction formula.

\section*{Acknowledgments}
This research did not receive any specific grant from funding agencies in the public, commercial, or not-for-profit sectors.


\begin{thebibliography}
\bibitem[\protect\citeauthoryear{Balazs {\it et~al}.}{2011}]{Balazs2011}
Balazs, P. and D\"orfler, M. and Holighaus, N. and Jaillet, F. and Velasco G~A. (2011). Theory, implementation and applications of nonstationary {Gabor} frames. {\it Journal of Computational and Applied Mathematics}, {\bf 236(6)}: 1481--1496. \url{10.1016/j.cam.2011.09.011}

\bibitem[\protect\citeauthoryear{Bhattacharyya {\it et~al}.}{2017}]{EEGEWT}
Bhattacharyya, A. and Gupta, V. and Pachori, R.~B. (2017). Automated identification of epileptic seizure EEG signals using empirical wavelet transform based {Hilbert} marginal spectrum. {\it Proceedings of the 22nd International Conference on Digital Signal Processing (DSP)}. \url{10.1109/icdsp.2017.8096122}

\bibitem[\protect\citeauthoryear{Boashash}{1992}]{Boashash1992}
Boashash, B. (1992). Estimating and interpreting the instantaneous frequency of a signal - part 1: Fundamentals. {\it Proceedings of the IEEE}, {\bf 80(4)}:520--538. \url{10.1109/5.135376}

\bibitem[\protect\citeauthoryear{Boashash}{1992}]{Boashash1992a}
Boashash, B. (1992). Estimating and interpreting the instantaneous frequency of a signal - part 2: Algorithms and applications. {\it Proceedings of the IEEE}, {\bf 80(4)}:540--568. \url{10.1109/5.135378}

\bibitem[\protect\citeauthoryear{Cattani}{2008}]{Shannon1}
Cattani, C. (2008). Shannon wavelets theory. {\it Mathematical Problems in Engineering}, {\bf Article ID 164808}:24 pages. \url{10.1155/2008/164808}

\bibitem[\protect\citeauthoryear{Chai and Shen}{2007}]{Chai2007}
Chai, A. and Shen, Z. (2007). Deconvolution: a wavelet frame approach. {\it Numerische Mathematik}, {\bf 106(4)}:529--587. \url{10.1007/s00211-007-0075-0}

\bibitem[\protect\citeauthoryear{Chambolle {\it et~al}.}{1998}]{Chambolle1998}
Chambolle, A., DeVore, R.~A. and Lee, N.~Y. and Lucier, B.~J. (1998). Nonlinear wavelet image processing: Variational problems, compression, and noise removal through wavelet shrinkage. {\it IEEE Transactions on Image Processing}, {\bf 7(3)}:319--335. \url{10.1109/83.661182}

\bibitem[\protect\citeauthoryear{Christensen}{2001}]{Christensen2001}
Christensen, O. (2001). Frames, {Riesz} bases, and discrete {Gabor}/wavelet expansions. {\it Bulletin of the American Mathematical Society}, {\bf 38(3)}:273--291. \url{10.1090/S0273-0979-01-00903-X}

\bibitem[\protect\citeauthoryear{Christensen}{2010}]{Christensen2010}
Christensen, O. and Laugesen, R.~S. (2010). Approximately dual frame pairs in {Hilbert} spaces and applications to {Gabor} frames. {\it Sampling Theory in Signal and Image Processing}, {\bf 9(1--3)}:77--89.

\bibitem[\protect\citeauthoryear{Daubechies}{1990}]{Daubechies1990}
Daubechies, I. (1990). The wavelet transform, time-frequency localization and signal analysis. {\it IEEE Transactions on Information Theory}, {\bf 36(5)}:961--1005. \url{10.1109/18.57199}

\bibitem[\protect\citeauthoryear{Daubechies}{1992}]{Daubechies1992}
Daubechies, I. (1992). Ten lectures on wavelets. {\it Society for Industrial and Applied Mathematics}, Philadelphia, PA, USA.

\bibitem[\protect\citeauthoryear{Daubechies {\it et~al}}{2011}]{Daubechies2011}
Daubechies, I. and Lu, J. and Wu, H-T. (2011). Synchrosqueezed wavelet transforms: An empirical mode decomposition-like tool. {\it Applied and Computational Harmonic Analysis}, {\bf 30(2)}:243--261. \url{10.1016/j.acha.2010.08.002}

\bibitem[\protect\citeauthoryear{Donoho and Stark}{1989}]{Donoho1989}
Donoho D.~L. and Stark, P.~B. (1989). Uncertainty principles and signal recovery. {\it SIAM Journal on Applied Mathematics}, {\bf 49(3)}:906--931.
\url{10.1137/0149053}

\bibitem[\protect\citeauthoryear{Dragomiretskiy and Zosso}{2014}]{Dragomiretskiy2014}
Dragomiretskiy, K. and Zosso, D. (2014). Variational mode decomposition. {\it IEEE Transactions on Signal Processing}, {\bf 62(3)}:531--544. \url{10.1109/TSP.2013.2288675}

\bibitem[\protect\citeauthoryear{Flandrin and Gon\c{c}alv\`es}{2004}]{Flandrin2004a}
Flandrin, P. and Gon\c{c}alv\`es, P. (2004). Empirical mode decompositions as data-driven wavelet-like expansions. {\it International Journal of Wavelets, Multiresolution and Information Processing}, {\bf 2(4)}:477--496. \url{10.1142/S0219691304000561}

\bibitem[\protect\citeauthoryear{Gilles}{2013}]{Gilles2013}
Gilles, J. (2013). Empirical wavelet transform. {\it IEEE Trans on Signal Processing}, {\bf 61(16)}:3999--4010. \url{10.1109/TSP.2013.2265222}

\bibitem[\protect\citeauthoryear{Gilles and Heal}{2014}]{Gilles2014a}
Gilles, J. and Heal, K. (2014). A parameterless scale-space approach to find meaningful modes in histograms - application to image and spectrum segmentation. {\it International Journal of Wavelets, Multiresolution and Information Processing}, {\bf 12(6)}:1450044--1--1450044--17. \url{10.1142/S0219691314500441}

\bibitem[\protect\citeauthoryear{Gilles and Osher}{2016}]{Gilles2016a}
Gilles, J. and Osher, S. (2016). Wavelet burst accumulation for turbulence mitigation. {\it Journal of Electronic Imaging}, {\bf 25(3)}:033003--1--033003--9. \url{10.1117/1.JEI.25.3.033003}

\bibitem[\protect\citeauthoryear{Gilles {\it et~al}}{2014}]{Gilles2013a}
Gilles, J. and Tran, G. and Osher, S. (2014). {2D} Empirical transforms. Wavelets, Ridgelets and Curvelets Revisited. {\it SIAM Journal on Imaging Sciences}, {\bf 7(1)}:157--186. \url{10.1137/130923774}

\bibitem[\protect\citeauthoryear{Hao and Jun-Hai}{2013}]{lfm}
Hao, C. and Jun-Hai, G. (2013). The application of data-driven {TF} analysis methods in {LFM} signal parameter estimation. {\it Proceedings of TENCON 2013 Conference}, 1--4, October. \url{10.1109/TENCON.2013.6718885}

\bibitem[\protect\citeauthoryear{Gilles}{2011}]{Hou2011}
Hou, T.~Y. and Shi, Z. (2011). Adaptive data analysis via sparse time-frequency representation. {\it Advances in Adaptive Data Analysis}, {\bf 2}:1--28. \url{10.1142/S1793536911000647}

\bibitem[\protect\citeauthoryear{Hramov {\it et~al}}{2015}]{Hramov2015}
Hramov, A.~E. and Koronovskii, A.~A. and Makarov, V.~A. and Pavlov, A.~N. and Sitnikova, E. (2015). {\it Wavelets in neuroscience}. Springer series in synergetics. Springer. \url{10.1007/978-3-662-43850-3}

\bibitem[\protect\citeauthoryear{Hua {\it et~al}}{2015}]{wind}
Hua, J. and Wangb, J. and Mac, K. (2015). A hybrid technique for short-term wind speed prediction. {\it Energy}, {\bf 81}:563--574. \url{10.1016/j.energy.2014.12.074}

\bibitem[\protect\citeauthoryear{Huang {\it et~al}}{1998}]{Huang1998}
Huang, N.~E. Shen, Z. and Long, S.~R. and Wu, M.~C. and Shih, H.~H. and Zheng, Q. and Yen, N.~C. and Tung, C.~C. and Liu, H.~H. (1998). The empirical mode decomposition and the {Hilbert} spectrum for nonlinear and non-stationary time series analysis. {\it Proc. Royal Society London A.}, {\bf 454}:903--995. \url{10.1098/rspa.1998.0193}

\bibitem[\protect\citeauthoryear{Huang {\it et~al}}{2018}]{Huang2018}
Huang, Y. and De~Bortoli, V. and Zhou, F. and Gilles, J. (2018). Review of wavelet-based unsupervised texture segmentation, advantage of adaptive wavelets. {\it IET Image Processing Journal}, {\bf 12(9)}:1626--1638. \url{10.1049/iet-ipr.2017.1005}

\bibitem[\protect\citeauthoryear{Huang {\it et~al}}{2019}]{Huang2019}
Huang, Y. and Zhou, F. and Gilles, J. (2019). Empirical curvelet based fully convolutional network for supervised texture image segmentation. {\it Neurocomputing}, {\bf 349}:31--43. \url{10.1016/j.neucom.2019.04.021}

\bibitem[\protect\citeauthoryear{Jaffard {\it et~al}}{2001}]{Jaffard2001}
Jaffard, S. and Meyer, Y. and Ryan, R.~D. (2001). Wavelets: Tools for Science and Technology. {\it SIAM}.

\bibitem[\protect\citeauthoryear{Kumara and Saini}{2014}]{ecg}
Kumara, R. and Saini, I. (2014). Empirical Wavelet Transform based {ECG} signal compression. {\it IETE Journal of Research}, {\bf 60(6)}:423--431.
\url{10.1080/03772063.2014.963173}

\bibitem[\protect\citeauthoryear{Li {\it et~al}}{2014}]{speech}
Li, Y. and Xue, B. and Hong, H. and Zhu, X. (2014). Instantaneous pitch estimation based on {Empirical Wavelet Transform}. {\it Proceedings of the 19th IEEE International Conference on Digital Signal Processing}, 250--253. \url{10.1109/ICDSP.2014.6900838}

\bibitem[\protect\citeauthoryear{Liu {\it et~al}}{2016}]{seismic}
Liu, W. and Cao, S. and Chen. Y. (2016). Seismic time-frequency analysis via empirical wavelet transform. {\it IEEE Geoscience and Remote Sensing Letters}, {\bf 13(1)}:28--32. \url{10.1109/lgrs.2015.2493198}

\bibitem[\protect\citeauthoryear{Livens {\it et~al}}{1997}]{Livens1997}
Livens, S. and Scheunders, P. and Van~de~Wouwer, G. and Van Dyck, D. (1997). Wavelets for texture analysis, an overview. {\it Image Processing and Its Applications, 1997., IET Sixth International Conference on}, {\bf volume~2}:581--585, Dublin. \url{10.1049/cp:19970958}

\bibitem[\protect\citeauthoryear{Mallat}{2009}]{Mallat2009}
Mallat, S. (2009). A Wavelet Tour of Signal Processing - A sparse way. {\it Academic Press}, Elsevier, third edition edition.

\bibitem[\protect\citeauthoryear{Morizet and Gilles}{2008}]{Gilles2008a}
Morizet, N. and Gilles, J. (2008). A new adaptive combination approach to score level fusion for face and iris biometrics combining wavelets and statistical moments. {\it International Symposium on Visual Computing (ISVC)}, Las Vegas, USA, December. \url{10.1007/978-3-540-89646-3_65}

\bibitem[\protect\citeauthoryear{Ocak}{2009}]{Ocak2009}
Ocak, H. (2009). Automatic detection of epileptic seizures in EEG using discrete wavelet transform and approximate entropy. {\it Expert Systems with Applications}, {\bf 36}:2027--2036. \url{10.1016/j.eswa.2007.12.065}

\bibitem[\protect\citeauthoryear{Ricaud {\it et~al}}{2013}]{Ricaud2013}
Ricaud, B. and Stempfel, G. and Torr\'esani, B. and Wiesmeyr, C. and Lachambre, H. and Onchis, D. (2013). An optimally concentrated {Gabor} transform for localized time-frequency components. {\it Advances in Computational Mathematics}. \url{10.1007/s10444-013-9337-9}

\bibitem[\protect\citeauthoryear{Shen}{2010}]{Shen2010}
Shen, Z. (2010). Wavelet frames and image restorations. In Rajendra~Bhatia eds, editor, {\em Proceedings of the International
  Congress of Mathematicians}, {\bf 4}:2834--2863, Hyderabad, India, Hindustan Book Agency.

\bibitem[\protect\citeauthoryear{Starck and Bijaoui}{1994}]{Startck1994}
Starck J.~L. and Bijaoui, A. (1994). Filtering and deconvolution by the wavelet transform. {\em Signal Processing}, {\bf 35}:195--211. \url{10.1016/0165-1684(94)90211-9}

\bibitem[\protect\citeauthoryear{Tang}{2009}]{Tang2009a}
Tang, Y.~Y. (2009). Wavelet theory approach to pattern recognition, {\it Series in Machine Perception and Artificial Intelligence}, {\bf 74}, World Scientific Publishing. \url{10.1142/7324}

\bibitem[\protect\citeauthoryear{Thirumala {\it et~al}}{2015}]{singlephase}
Thirumala, K. and Umarikar, A.~C. and Jain, T. (2015). Estimation of single-phase and three-phase power-quality indices using empirical wavelet transform. {\it IEEE Transactions on power delivery}, {\bf 30(1)}. \url{10.1109/TPWRD.2014.2355296}

\bibitem[\protect\citeauthoryear{Triebel}{2000}]{Triebel2}
Triebel H. (2000). Theory of Function Spaces II. {\it Monographs in Mathematics}. Birkha\"user Verlag.

\bibitem[\protect\citeauthoryear{Triebel}{2006}]{Triebel3}
Triebel, H. (2006). Theory of Function Spaces III. {\it Monographs in Mathematics}. Birkha\"user Verlag.

\bibitem[\protect\citeauthoryear{Vedel}{2009}]{Vedel2009}
Vedel, B. Flat wavelet bases adapted to the homogeneous {Sobolev} spaces. {\it Mathematische Nachrichten}, {\bf 282(1)}:104--124. \url{10.1002/mana.200810725}

\end{thebibliography}
\end{document}